\documentclass[12pt,leqno]{article}
\usepackage{amsmath,amssymb,theorem,epic} \title{Decomposition of
  tensor products \\ of modular irreducibles for $\SL_2$}

\author{Stephen Doty and Anne Henke}

\newenvironment{pf}{\noindent{\em Proof.}}{\hfill$\square$\par\medskip} 
\makeatletter 
\renewcommand{\subsection}{\@startsection{subsection}{2}{0mm}{\baselineskip}{-2\fontdimen2\font}{\normalfont\normalsize\bfseries}}
\makeatother
 
\theoremstyle{plain} 
\newtheorem{thm}[equation]{Theorem}
\newtheorem{cor}[equation]{Corollary}
\newtheorem{lem}[equation]{Lemma}

\newtheorem{conj[equation]}{Conjecture}
{\theorembodyfont{\rmfamily}}
{\theorembodyfont{\rmfamily}}
{\theorembodyfont{\rmfamily}} 

\numberwithin{equation}{section}

\newcommand{\Z}{{\mathbb Z}}

\newcommand{\Ext}{\operatorname{Ext}}

\newcommand{\ind}{\operatorname{ind}}

\newcommand{\soc}{\operatorname{soc}}

\renewcommand{\ker}{\operatorname{Ker}}
\newcommand{\GL}{{\sf GL}}
\newcommand{\SL}{{\sf SL}}
\newcommand{\boldu}{{\mathbf u}}

\begin{document}
\maketitle
\begin{abstract}
  We use tilting modules to study the structure of the tensor product
  of two simple modules for the algebraic group $\SL_2$, in positive
  characteristic, obtaining a twisted tensor product theorem for its
  indecomposable direct summands. Various other related results are
  obtained, and numerous examples are computed.
\end{abstract}

\section*{Introduction}

We study the structure of $L\otimes L^\prime$ where $L,L^\prime$ are
simple modules for the algebraic group $\SL_2 = \SL_2(k)$ over an
algebraically closed field $k$ of positive characteristic $p$. The
solution to this problem is well-known and easily obtained in
characteristic zero, but in positive characteristic the problem 
is significant. 

Given $L \otimes L'$, the initial question is to describe its
indecomposable direct summands.  This is answered in Theorem 2.1. It
turns out that each such direct summand is expressible as a twisted
tensor product of certain ``small'' indecomposable tilting modules
where the structure of the latter is completely understood (see Lemma
\ref{lem:Fundamental}). We note that for $p=2$ the module $L\otimes
L'$ is always indecomposable, in contrast to what happens for $p$ odd.

The indecomposable summands themselves are always contravariantly
self-dual, with simple socle (and head), and they occur as
subquotients of tilting modules (see Theorem \ref{tilt:subquotients}).
On the other hand, each tilting module occurs as a direct summand of
some $L\otimes L'$ (see Theorem \ref{tilt:summand}).  These results
provide the starting point for calculating the examples in Section 5
and 6.

In fact, we suggest that the reader start by browsing through the
examples in Sections 5 and 6.  The close relationship between tensor
products of simples and tilting modules will be apparent from these
examples.  Since examples of tilting module structure are rare in the
literature, these computations should be of independent interest.

Our general results are given in Sections 2 through 4.  In particular,
we classify precisely which indecomposable summands of $L\otimes L'$
are tilting, and we obtain a result expressing certain indecomposable
tilting modules as tensor products of two simple modules (usually in
more than one way). 

In Sections \ref{sec:tpstructure} and 4 we study in detail the tensor
products $L \otimes L(1)$ where $L$ is arbitrary.  Here we obtain a
class of uniserial and biserial tilting modules and their ``shifts''.
The methods used here can be applied to obtain the structure of all
tensor products $L \otimes L'$ where $L'=L(a)$ with $a \leq p-1$,
although we formulate the precise statement only for $L \otimes L(2)$.

The results of the paper are related to results of Alperin
\cite{Alperin}, Brundan and Kleshchev \cite{BK:translation}, and
Erdmann and Henke \cite{EH:RingelDuality, EH:uniserial}. Our main
technical tools are Steinberg's tensor product theorem and Donkin's
interpretation \cite{Donkin:tilt}, in the context of algebraic groups,
of Ringel's theory \cite{R:tilting} of tilting modules.

\section{Preliminaries} \label{Prelim}
The set $X=X(T)$ of weights for a maximal torus $T$ in the algebraic
group $\SL_2$ will be identified with the set $\Z$ of integers, as
usual. Then dominant weights correspond to nonnegative integers. If
$r$ is such then we write $L(r)$ for the simple $\SL_2$-module of
highest weight $r$ and $\Delta(r)$ for the Weyl module of that same
highest weight.  We write $\nabla(r)$ for the transpose
(contravariant) dual of $\Delta(r)$.  By a theorem of Ringel (see
\cite{R:tilting}), there exists a unique indecomposable module
$T(r)$ of highest weight $r$ such that $T(r)$ has both
$\Delta$-filtration and $\nabla$-filtration. The modules $T(r)$ are
the (partial) tilting modules.

If a module $M$ has a composition series 
$
0 = M_0 \le M_1 \le \cdots \le M_k = M
$
with simple factors $S_i \cong M_i/M_{i-1}$ for $i=1, \dots, k$ then
we denote that composition series by writing $$[S_1, S_2, \dots, S_k].$$

We begin with some easy lemmas on tilting modules.  Our first result
describes the module structure of certain small tilting modules, which
will turn out to be the basic material out of which all tilting
modules and tensor products of simples are built up.

\begin{lem}\label{lem:Fundamental}
\begin{itemize}
\item[(a)] For $0 \le u \le p-1$ we have $T(u)=L(u)=\nabla(u)=\Delta(u)$.  
\item[(b)] For $p \le u \le 2p-2$ the module $T(u)$ is uniserial
and its unique composition series has the form $[L(2p-2-u), L(u),
L(2p-2-u)]$. Moreover, $T(u)$ is a non-split extension of
$\Delta(2p-2-u)$ by $\Delta(u)$ (or, dually, of $\nabla(u)$ by
$\nabla(2p-2-u)$).
\end{itemize}
\end{lem}
\begin{pf}
  Recall that $\nabla(r) \cong S^r(E)$, the $r$th symmetric power of
  the natural module $E$.  Part (a) follows immediately from the fact
  that $\nabla(u) = \Delta(u) = L(u)$ is simple for $0 \le u \le p-1$.
  This fact is well-known and follows for instance from the strong
  linkage principle, or from known results \cite{Doty:thesis} on the
  structure of $S^r(E)$.

Part (b) is a special case of \cite[Proposition 2.3]{EH:uniserial}.
Or: Write $G=\SL_2$ and set $G_1 = \ker F$ (the Frobenius kernel).  As
observed by Donkin \cite[\S2, Example 1]{Donkin:tilt}, for $0\le
\lambda \le p-1$ we have an isomorphism $T(2p-2-\lambda) \cong
Q(\lambda)$, where $Q(\lambda)$ is the unique (up to $G$-isomorphism)
$G$-module such that $Q(\lambda)|_{G_1}$ is isomorphic with the
projective cover of $L(\lambda)|_{G_1}$.  The existence of this
$G$-module lift follows from results of Jantzen, extending earlier
results of Ballard.  Now we can apply \cite[II, 11.4 Prop.]{Jantzen}
to compute the formal character of $Q(\lambda)$. Setting
$u=2p-2-\lambda$ and restricting to $0\le \lambda\le p-2$, we can
apply the aformentioned description of the $\nabla$'s as symmetric
powers to obtain part (b) of the lemma.
\end{pf}

Let us call the tilting modules described in the preceding result {\em
fundamental}. As we shall see (in Lemma \ref{lem:TwistedTilting}
ahead) any tilting module for $\SL_2$ can be expressed as a twisted
tensor product of fundamental tilting modules.  Moreover, in Theorem
\ref{thm:Tensor} we shall see that the indecomposable summands of $L
\otimes L'$ can also be expressed as such a twisted tensor
product. (The indecomposable summands are not necessarily tilting,
however.)  For $0\le u \le 2p-2$, we denote by $\widetilde{u}$ the
highest weight of the socle (and head) of $T(u)$, so that
\begin{eqnarray}  \label{eqn1}
\widetilde{u}&=&
\begin{cases}
u & \text{if $u \le p-1$,}\\
2p-2-u & \text{otherwise}.
\end{cases}
\end{eqnarray}

Most cases of the next lemma appear already in \cite[Lemma 4]
{EH:RingelDuality}.

\begin{lem}\label{lem:SmallTensor}
Let $L,L'$ be two simple modules in the bottom alcove, i.e.\ their
highest weights are inclusively between $0$ and $p-1$. Then $L\otimes
L'$ is tilting, and isomorphic with the direct sum of $T(u)$ as $u$
varies over a set $W(L,L')$ of weights which can be computed as
follows. Let $r$ (resp., $s$) be the larger (resp., smaller) of the
highest weights of $L,L'$.  List the weights $r+s$, $r+s - 2, \dots,
r-s$. For each $u \ge p$ on this list, strike out the number $2p-2-u$
from the list. What remains is the set $W(L,L')$.  In other words, if
$S = \{r+s-2i\}_{i=0}^{s}$ then
$$
W(L,L')= S - \{2p-2-u \mid u \in S, u\ge p \}.
$$
In particular, $L \otimes L'$ is indecomposable if and only if
$s=0$ or $(r,s)=(p-1,1)$.
\end{lem}
\begin{pf}
  The module $L(r) \otimes L(s)$ is tilting. It has a
  $\nabla$-filtration with sections $\nabla(r+s)$, $\nabla(r+s-2),
  \dots, \nabla(r-s)$.  Each of these sections has either one or two
  composition factors.  The statement about $W(L,L')$ now follows from
  Lemma \ref{lem:Fundamental} and some simple bookkeeping.

  It is clear that $L\otimes L'$ is indecomposable in case $s=0$ or
  $(r,s) = (p-1,1)$.  For the converse, note that if $s=1$ and $r\ne
  p-1$ then $W(L,L') = \{r+1, r-1\}$. If $s>1$ then $S$ contains at
  least three elements. If $r+s-2 \notin W(L,L')$ then $r+s=p$ and
  $r-s \in W(L,L')$, so $W(L,L')$ either contains both $r+s, r-s$ or
  both $r+s, r+s-2$. This proves the indecomposability claim.
\end{pf} 
For later use we set $W(a,b) = W(L,L')$ where $L=L(a)$, $L'=L(b)$ and
$0\le a,b \le p-1$.  Note that $W(a,b)=W(b,a)$.

%\begin{remark}
%In case $p=2$, all $L\otimes L'$ are indecomposable, for $L,L'$ in the
%bottom alcove.
%\end{remark}

We shall also need the following lemma from \cite[Lemma
5]{EH:RingelDuality}, which describes how the tilting modules are
built up as twisted tensor products of fundamental modules.  We use a
superscript $F^i$ applied to a module to indicate twisting by the
$i$th iterate of the Frobenius morphism.
\begin{lem}\label{lem:TwistedTilting}
Let $u \in \Z_{\ge0}$.  Then $u$ can be written uniquely in the form
$u = \sum_{i=0}^m u_ip^i$ where $p-1 \le u_i \le 2p-2$ for $i < m$ and
$0 \le u_m \le p-1$.  Then $T(u) \simeq \bigotimes_{i=0}^m
T(u_i)^{F^i}$.
\end{lem}

\section{The indecomposable summands of $L \otimes L'$}\label{Indec}
The next result follows easily from Lemmas \ref{lem:Fundamental} and
\ref{lem:SmallTensor}, combined with Steinberg's tensor product
theorem.  For $r\in \Z_{\ge0}$, write $r=\sum_{i\ge0} \delta_i(r) p^i$
where $\delta_i(r) \in \{0,1, \dots, p-1\}$ (the $p$-adic expansion of
$r$). In the following, when writing $W(\delta_i(r), \delta_i(r'))$ we
mean the set described in Lemma \ref{lem:SmallTensor}. 

\begin{thm}\label{thm:Tensor}
Let $r, r'$ be arbitrary nonnegative integers.  The tensor product
$L(r)\otimes L(r')$ can be expressed as a direct sum of twisted tensor
products of fundamental tilting modules.  In fact, we have
$$
L(r) \otimes L(r') \simeq \bigoplus_\boldu \left( \bigotimes_i
T(u_i)^{F^i} \right)
$$
where $\boldu=(u_0, \dots, u_m)$ ranges over all elements of
the finite Cartesian product
$$
W=W(\delta_0(r), \delta_0(r')) \times W(\delta_1(r), \delta_1(r')) \times \cdots \times
W(\delta_m(r), \delta_m(r'))
$$
of the sets described in Lemma \ref{lem:SmallTensor}, and where $m$
is the $p$-adic length of the largest of $r$, $r'$.  Given $\boldu$ as
above, the corresponding indecomposable direct summand $J(\boldu) =
\bigotimes_{i=0}^m T(u_i)^{F^i}$ is always contravariantly self-dual,
with simple socle and head isomorphic with $L(\sum_{i=0}^m
\widetilde{u_i} p^i)$ where $\widetilde{u_i}$ is defined as in
equation (\ref{eqn1}).
\end{thm}

\begin{pf}
  For a given integer $r\ge0$ we have its $p$-adic expansion
  $r=\sum_{i\ge0} \delta_i(r)p^r$ where each $\delta_i(r)$ satisfies
  $0 \le \delta_i(r) \le p-1$. By Steinberg's tensor product theorem
  $L(r)\simeq \bigotimes_i L(\delta_i(r))^{F^i}$. Thus (using Lemma
  \ref{lem:SmallTensor}) we have
\begin{align*}
  L(r)\otimes L(r') &\simeq \bigotimes_i [L(\delta_i(r))\otimes 
  L(\delta_i(r'))]^{F^i}\\
  &\simeq \bigotimes_i \bigoplus_{u\in W(\delta_i(r),\delta_i(r'))}
  T(u)^{F^i}
\end{align*}
and the stated isomorphism follows by interchanging products and sums.
The statements about $J(\boldu)$ follow from Lemma \ref{lem:Fundamental}
and Steinberg's tensor product theorem.
\end{pf}

Note that the multiplicities of the composition factors of each
$J(\boldu)$ are computable from Lemma \ref{lem:Fundamental} and
Steinberg's tensor product theorem. 
See the examples in the last two sections.

%For example, let $p=3$ then $L(8)
%\otimes L(2)$ has $W=W(2,2) \times W(2,0)=\{ (2,2), (4,2) \}$  and 
%hence 
%\[
%L(8) \otimes L(2) = T(2) \otimes T(2)^{F} \oplus T(4) \otimes T(2)^{F}
%=L(8) \oplus [0,4,0]  \otimes T(2)^{F}
%\]
%as $T(4)=[0,4,0]$ (fundamental). On the other hand, $L(6) \otimes
%L(4)=T(1) \otimes T(3)^{F}$ with $T(3)=[1,3,1]$ (fundamental).
%Hence $T(4) \otimes L(2)^{F}$ has composition factors twice 6, twice
%4 and once 10. (It is the uniserial module $[6,4,10,4,6]$.) 

We can prove the following criterion, classifying the summands of the
tensor product that are tilting and/or irreducible:
\begin{thm}\label{thm:TiltingSummand}
  For $\boldu$, $m$ as in the preceding theorem, the indecomposable
  summand $J(\boldu)$ is a tilting module if and only if $u_i \ge p-1$
  for $i=0, 1, \dots, m-1$, in which case it is isomorphic with the
  tilting module $T(\sum_{i=0}^m u_i p^i)$. Moreover, $J(\boldu)$ is
  irreducible if and only if $u_i \le p-1$ for all $i=0, 1, \dots, m$,
  in which case it is isomorphic with $L(\sum_{i=0}^m u_i p^i)$.
\end{thm}

\begin{pf} The second claim follows immediately from Steinberg's
tensor product theorem and Lemma \ref{lem:Fundamental}.  The ``if''
part of the first claim follows immediately from Lemma
\ref{lem:TwistedTilting}.  For the converse assertion, suppose that
$J(\boldu)$ is tilting but there exists some index $j<m$ for which
$u_j<p-1$. The highest weight of $J(\boldu)$ is $u'=\sum_{i=0}^m
u_ip^i$, so $J(\boldu)\cong T(u')$.  We express $u'$ as in Lemma
\ref{lem:TwistedTilting} and decompose $T(u')$ into a twisted tensor
product of tilting modules.  Using Lemma \ref{lem:Fundamental} we can
compute the highest weight of the socle of $T(u')$. Now a computation
shows that the highest weight of the socle of $J(\boldu)$ will differ
from the highest weight of the socle of $T(u')$, and thus we arrive at
a contradiction. The details are as follows:

Let $u' =\sum_{i=0}^t v_i p^i$ be an admissable decomposition as in
Lemma \ref{lem:TwistedTilting}. As $\soc J(\boldu)= \soc T(u')$ we
obtain together with Theorem \ref{thm:Tensor}
\[
\sum_{i=0}^m\tilde{u}_i =\sum_{i=0}^t \tilde{v}_ip^i.
\] 
As these are $p$-adic decompositions we obtain
$\tilde{v}_i=\tilde{u}_i$.  Certainly $t \leq m+1$ and if $t=m+1$ then
$\tilde{v}_t=0=v_t$. We assume now that $t=m$. (For $t <m$ a similar
calculation leads to a contradiction.) Define
\begin{eqnarray*}
I=\{ i \mid u_i <p-1, i \neq t\} \mbox{ \hspace{.5cm}and \hspace{.5cm}}  R=\{ i \mid u_i \geq p-1, i \neq t\}.
\end{eqnarray*}
So by assumption $I$ is not empty. We obtain
\begin{eqnarray*}
\sum_{i=0}^m v_i p^i =u' &=& \sum_{i=0}^m u_i p^i =\sum_{i \in I} u_i p^i + \sum_{i \in R} u_i p^i +u_mp^m \\
&=& \sum_{i \in I} (2p-2-v_i) p^i + \sum_{i \in R} v_i p^i +u_mp^m,  
\end{eqnarray*}
which is equivalent to 
\begin{eqnarray*}
2 \sum_{i \in I} v_i p^i +v_mp^m &=& 2 \sum_{i \in I} (p-1) p^i +u_mp^m.
\end{eqnarray*}
Either $v_m=u_m <p-1$ or $u_m=2p-2-v_m$. In both cases we obtain
$v_i=p-1$ and hence $u_i=v_i=p-1$. This contradicts the assumption
that $I$ is not empty.
\end{pf}

We also have the following criterion for the tensor product $L(r)
\otimes L(r')$ to be indecomposable:

\begin{thm}\label{thm:Indecomposable}
With $m$ as above, $L(r)\otimes L(r')$ is indecomposable if and only
if, for each $i = 0, \dots, m$, either 
\begin{itemize}
\item[(a)] one of $\delta_i(r)$ or $\delta_i(r')$ is $0$, or
\item[(b)] $\delta_i(r) + \delta_i(r') = p$ and one of $\delta_i(r)$ or
$\delta_i(r')$ is $1$.
\end{itemize}
\end{thm}

\begin{pf}
Combine Lemma \ref{lem:SmallTensor} with Theorem \ref{thm:Tensor}. 
\end{pf}

Note that for $p=2$ the conditions of the preceding theorem
necessarily hold for any $r$, $r'$, so $L(r)\otimes L(r')$ is always
indecomposable when $p=2$.  By combining the previous two results we
obtain:
\begin{thm}\label{thm:IndecomposableTilting} 
  With $m$ as above, necessary and sufficient conditions for the
  tensor product $L(r)\otimes L(r')$ to be isomorphic to an
  indecomposable tilting module are that:
\begin{itemize}
\item[(a)] $p-1 \le \delta_i(r) + \delta_i(r') \le p$, for 
all $i=0,\dots, m-1$, and
\item[(b)] $\delta_m(r)+\delta_m(r')\le p$, and 
\item[(c)] when $\delta_i(r)+\delta_i(r')=p$, one of $\delta_i(r)$ or
$\delta_i(r')$ is equal to $1$, and 
\item[(d)] when $\delta_i(r)+\delta_i(r')<p$, one of $\delta_i(r)$ or 
$\delta_i(r')$ is equal to $0$. 
\end{itemize}
When these conditions hold we have the isomorphism $L(r)\otimes L(r')
\simeq T(r+r')$.
\end{thm}

\begin{pf}
This follows immediately from Theorems \ref{thm:TiltingSummand}
and \ref{thm:Indecomposable}.
\end{pf}

We also have the following criterion for when a given tilting module
can be factored as a tensor product of two simple modules.
\begin{thm}\label{thm:TiltingFactorization}
  Necessary and sufficient conditions for $T(u)$ to be factorizable as
  a tensor product of two simple modules are that when $u$ is
  expressed in the form $u = \sum_{i=0}^m u_i p^i$ with $u_i$ as in
  Lemma \ref{lem:TwistedTilting}, we have $u_i = p$ or $p-1$ for all
  $i=0, \dots, m-1$.  In that case, if $p>2$, there are precisely
  $2^m$ or $2^{m-1}$ such factorizations, depending on whether $u_m>0$
  or $u_m=0$, respectively. For $p=2$ there are only $2^t$ such
  factorizations, where $t = |\{i \mid u_i =p-1=1 \}|$.
\end{thm}

\begin{pf}
This follows from Theorem \ref{thm:IndecomposableTilting}.
\end{pf}

Note that for $p=2$ every indecomposable tilting module has a
factorization of the above form.  This is not true for odd primes: for
$p=3$ there are no integers $r$, $r'$ for which the module $T(4)$ is
isomorphic with $L(r) \otimes L(r')$. However, $T(4)$ is a direct
summand of $L(2) \otimes L(2)$, and in fact we have the following
statement in general.

\begin{thm}\label{tilt:summand}
Every indecomposable tilting module occurs as a direct summand of some
product $L \otimes L'$ for appropriate simple modules $L,L'$.
\end{thm}

\begin{pf}
Given $T(u)$, express $u = \sum_{i=0}^m u_i p^i$ with $p-1 \le u_i \le
2p-2$ for $i<m$ and $0\le u_m \le p-1$.  For each $i$ we can find
integers $r_i, s_i$ in the range $0 \le r_i, s_i \le p-1$ such that
$u_i \in W(r_i,s_i)$.  (One can, for instance, simply choose $r_i$,
$s_i$ to satisfy $r_i + s_i = u_i$.)  Then by Theorem \ref{thm:Tensor}
and Lemma \ref{lem:TwistedTilting} it follows that $T(u)$ is a direct
summand of $L(\sum_i r_i p^i) \otimes L(\sum_i s_i p^i)$.
\end{pf}

\begin{thm} \label{tilt:subquotients}
  The indecomposable summands of a tensor product of two simple
  modules $L(r)$ and $L(r')$ are always obtainable as subquotients of
  tilting modules, more precisely, of the direct summands of $T(r) \otimes
  T(r')$. In particular, if $L(r) \otimes L(r')$ is indecomposable
  then it is a subquotient of $T(r) \otimes T(r')$.
\end{thm}

\begin{pf}
  Given simple modules $L(r)$ and $L(r')$. As tilting modules have a
  $\Delta$-filtration we can choose submodules $M(r)$ and $M(r')$ of
  $T(r)$ and $T(r')$ respectively such that $L(r)$ and $L(r')$ occur
  in the head of $M(r)$ and $M(r')$ respectively. Hence $L(r) \otimes
  L(r')$ is a quotient of $M(r) \otimes M(r')$. Moreover, $M(r)
  \otimes M(r')$ is a submodule of the tilting module $T(r) \otimes
  T(r')$. As the head of each indecomposable summand of $L(r) \otimes
  L(r')$ is simple (by Theorem~\ref{thm:Tensor}, we obtain each such
  summand as a subquotient of a tilting module $T$ where $T$ is a
  direct summand of the tilting module $T(r) \otimes T(r')$. The last
  claim follows by comparing heighest weights.
\end{pf}

%%%%%%%%%%%%%%%%%%%%%%%%%%%%%%%%%%%%%%%%%%%%%%%%%%%%%%%%%%%%%%%%%%%%%%%
\section{Structure of tensor products}\label{sec:tpstructure}
In this section we study the stucture of the tensor product of two
simple modules for $\SL_2$. This is a more difficult problem than those
considered so far, and our results are less than definitive. 

Let $r\ge 0$. We use the term {\it biserial} to refer to a
module that has precisely two distinct composition series. Our next
two results describe the structure of $M=L(r)\otimes L(1)$ for
arbitrary $r$.  These results are related to results of Brundan and
Kleshchev \cite{BK:translation}. The proof of the next two
theorems is given in the next section.
\begin{thm}\label{thm:LtensorL1} 
  The following cases occur for $M=L(r)\otimes L(1)$:
  \begin{itemize} 
  \item[(a)] If $r \equiv 0 \pmod{p}$ then $M$ is simple and
    isomorphic with $L(r+1)$.
  \item[(b)] If $r \equiv -1 \pmod{p}$ then $M$ is an indecomposable
    module.  Moreover, there exist unique positive integers $t$, $a$
    such that $r \equiv ap^t-1 \pmod{p^{t+1}}$ for $1 \leq a<p,$ and
    $M$ is uniserial when $a=1$ and is biserial otherwise.
  \item[(c)] In all other cases (that is $r \not\equiv 0,-1 \pmod{p}$)
    $M$ is isomorphic to the direct sum $L(r+1) \oplus L(r-1)$.
\end{itemize}
\end{thm}

The proof for the cases (a) and (c) follows from Theorem
\ref{thm:Tensor}, Lemma \ref{lem:Fundamental} and Steinberg's tensor
product theorem.  Note that case (c) cannot occur when $p=2$ and also
the biserial part of case (b) never occurs when $p=2$.  Our next
result describes the precise structure of $M=L(r)\otimes L(1)$ in case
(b).
\begin{thm}\label{thm:Structure} 
  Suppose that $r \equiv -1 \pmod{p}$.  Then $r=(ap^t-1)+p^{t+1}k$ for
  $1 \leq a \leq p-1$ and $k\in\Z$. Set $r_0 = r+1$ and $r_i =
  r+1 - 2p^{i-1}$ for $i \in \Z_{> 0}$. Then $M=L(r)\otimes L(1)$ has
  the following precise structure:
\begin{itemize}
\item[(a)] In case $a=1$ the module $M$ is the uniserial module of the
form
$$
[L(r_1), \dots, L(r_t), L(r_0), L(r_t), \dots, L(r_1)].
$$
\item[(b)] In all other cases $M$ is the biserial module 
with the structure diagram:
$$M = 
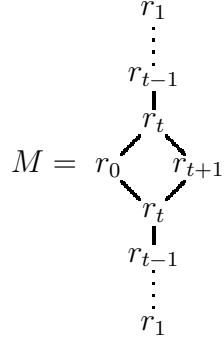
\begin{minipage}{12mm}
  \setlength{\unitlength}{0.6mm}
\begin{picture}(20,40)(-15,10)
\thicklines
\drawline[40](3,3)(7,7)
\drawline[40](-3,3)(-7,7)

\drawline[40](-7,13)(-3,17)
\drawline[40](7,13)(3,17)

\drawline[40](0,23)(0,27)
\drawline[40](0,-3)(0,-7)
\dottedline{3}(0,33)(0,41)
\dottedline{3}(0,-13)(0,-21)

\put(0,-10){\makebox(0,0){$r_{t-1}$}}
\put(0,-25){\makebox(0,0){$r_{1}$}}

\put(0,0){\makebox(0,0){$r_t$}}
\put(-10,10){\makebox(0,0){$r_0$}}
\put(10,10){\makebox(0,0){$r_{t+1}$}}
\put(0,20){\makebox(0,0){$r_t$}}

\put(0,30){\makebox(0,0){$r_{t-1}$}}
\put(0,45){\makebox(0,0){$r_{1}$}}

\end{picture}
\vspace{18mm}
\end{minipage}
$$
\end{itemize}
Finally, when $k=0$ in both cases $M$ is isomorphic with the
indecomposable tilting module of highest weight $r_0$ and has a
$\nabla$-filtration the factors of which are the uniserial modules
$\nabla(r_0)$ and $\nabla(r_1)$.
\end{thm}
(The structure diagram of $M$ follows the conventions of Alperin
\cite{Alperin:diagrams}, in which vertices correspond to simple
factors and edges to non-split extensions.  In the diagram we identify
simple composition factors $L(r_i)$ with their highest weight $r_i$.)

To describe the tensor products $M=L(r)\otimes L(a)$ for $2 \leq a
\leq p-1$ and $p >a$ is an application of the preceding results for
$L(r)\otimes L(1)$. We illustrate this on the tensor products $M=L(r)
\otimes L(2)$.  Set $E = L(1)$, the natural module for $\SL_2$.

Using the associativity of tensor products $(L \otimes E) \otimes E
\simeq L \otimes (E \otimes E)$, the module $L \otimes E \otimes E$
can be calculated in two ways. As $E \otimes E \simeq L(0) \oplus
L(2)$ for $p>2$, this provides a way to calculate $M=L(r) \otimes
L(2)$ for $p >2$. The result (following immediately from the preceding
results) is the following:
\begin{thm}\label{thm:LtensorL2}
  For $p=3$ the following cases occur for $M=L(r) \otimes L(2)$:
\begin{itemize}
\item[(a)] If $r \equiv 0 \pmod{p}$ then $M$ is isomorphic to $L(r+2)$.
\item[(b)] If $r \equiv 1 \pmod{p}$ then $M$ is isomorphic to $L(r+1)
  \otimes L(1) \simeq T(ap^t) \otimes L(k)^{F^{t+1}}$ where
  $r+1=(ap^t-1) +kp^{t+1}$ and $a=1,2$.
\item[(c)] If $r \equiv -1 \pmod{p}$ then $M$ is isomorphic to
  $T(ap^t+1) \otimes L(k)^{F^{t+1}} \oplus L(r)$ where $r=(ap^t-1)
  +kp^{t+1}$ and $a=1,2$.
\end{itemize}
\end{thm}
\begin{thm}\label{thm:L2Structure}
For $p \geq 5$ the following cases occur for $M=L(r) \otimes L(2)$:
\begin{itemize}
\item[(a)] If $r \equiv 0 \pmod{p}$ then $M$ is isomorphic to $L(r+2)$.
\item[(b)] If $r \equiv 1 \pmod{p}$ then $M$ is isomorphic to $L(r+2)
  \oplus L(r)$.
\item[(c)] If $r \equiv -1 \pmod{p}$ then $M$ is isomorphic to
  $T(ap^t+1) \otimes L(k)^{F^{t+1}} \oplus L(r)$ where $r=(ap^t-1)
  +kp^{t+1}$ and $1 \leq a \leq p-1$.
\item[(d)] If $r \equiv -2 \pmod{p}$ then $M$ is isomorphic to
  $T(ap^t) \otimes L(k)^{F^{t+1}} \oplus L(r-2)$ where $r+1=(ap^t-1)
  +kp^{t+1}$ and $1 \leq a \leq p-1$.
\item[(e)] In all other cases $M$ is isomorphic to $L(r+2) \oplus L(r)
  \oplus L(r-2)$.
\end{itemize}
\end{thm}
It should be noted that the submodule lattices of these tensor products
can be deduced from Theorem \ref{thm:Structure}.
\section{Proof of Theorems \ref{thm:LtensorL1} and
\ref{thm:Structure}} \label{Proof}
For a positive integer $t$ let $X_t$ be the set of $p^t$-restricted
weights for $\SL_2$; that is, the set of all $r \in \Z$ such that $0\le
r < p^t$. 

\begin{lem}\label{lem:1}
  Let $M$ be an $\SL_2$-module such that the highest weight of all
  composition factors of $M$ belongs to $X_t$.  Then for any $\mu\in
  \Z_{\ge0}$ the functor $\mathrm{(-)}\otimes L(\mu)^{F^{t+1}}$ induces
  an isomorphism of submodule lattices between $M$ and $M\otimes
  L(\mu)^{F^{t+1}}$. Moreover, the composition factors of $M\otimes
  L(\mu)^{F^{t+1}}$ have the form $L(\omega+p^{t+1}\mu)$ as $\omega$
  ranges over the set of highest weights of composition factors of
  $M$.
\end{lem}

\begin{pf}
This follows immediately from Steinberg's tensor product theorem,
using induction on the length of a composition series for $M$.
\end{pf}

Recall the notation $\delta_i(m)$ for the $i$th digit in the $p$-adic
expansion of $m$, so that $m=\sum_{i\ge0} \delta_i(m)p^i$ with $0\le
\delta_i(m) < p$ for all $i$.  We have the following elementary
result:
\begin{lem}\label{lem:2}
Let $r>0$ be an integer with $r\equiv -1 \pmod{p}$. There exist unique
positive integers $t$, $a$ such that $a<p$ and $r\equiv ap^t-1
\pmod{p^{t+1}}$.
\end{lem}

\begin{pf}
Let $t$ be the smallest index $i$ such that $\delta_i(r)$ is different
from $p-1$. By the hypothesis $t>0$. Since $\delta_i(r)=p-1$ for
$i=0,\dots, t-1$ we have
\begin{align*}
r & \equiv  \sum_{i=0}^{t-1} (p-1)p^i + \delta_t(r)p^t  \pmod{p^{t+1}} \\
  & \equiv  (p^t-1) + \delta_t(r)p^t  \pmod{p^{t+1}} \\
  & \equiv  ap^t-1  \pmod{p^{t+1}} 
\end{align*}  
where $a = 1 + \delta_t(r)$. This proves the existence of $t$, $a$ with
the stated properties.

If $t_1$, $t_2$, $a_1$, $a_2$ are positive integers with $a_i<p$ such
that $r\equiv a_ip^{t_i}-1 \pmod{p^{t_i+1}}$ for $i=1, 2$ then 
$r=(a_i p^{t_i}-1)+k_ip^{t_i+1}$ for  integers $k_i$. Hence
$$
a_1 p^{t_1} + k_1p^{t_1+1} = a_2 p^{t_2} + k_2p^{t_2+1}.
$$
If $t_1<t_2$ then upon dividing by $p^{t_1}$ we obtain the
contradiction that $p$ divides $a_1$. Similarly if $t_1>t_2$. Thus
$t_1 = t_2$ and so
$
a_1 + k_1 p = a_2 + k_2 p
$ 
from which we obtain immediately that $a_1=a_2$, $k_1=k_2$. This
proves the uniqueness statement.
\end{pf}

\begin{lem}\label{lem:3}
  For $r \in \Z_{\ge0}$ of the form $r=ap^t-1$ for $a \in \{1,\dots,
  p-1\}$ we have $\nabla(r) = L(r) = \Delta(r)$.
\end{lem}

\begin{pf}
  First, note that $\nabla(p^t-1,0) = L(p^t-1,0) = \Delta(p^t-1,0)$.
  This is well-known. It follows for $p$ odd from \cite[II,
  3.19(4)]{Jantzen}.  For general $p$ it follows for instance from
  \cite[Theorem 2]{Doty:linkage}.

  Moreover, from \cite[Theorem 1]{Doty:linkage} it follows that
  $\nabla_t(p^t-1,0) = L_t(p^t-1,0)$ where $\nabla_t(\lambda)$ (resp.,
  $L_t(\lambda)$) is the induced $G_tB$-module $\ind_B^{G_tB} \lambda$
  (resp., the simple $G_tB$-module of highest weight $\lambda$). Here
  $G=\SL_2$ and $G_t$ (resp., $G_tB$) is the kernel (resp., the
  inverse image of $B$) under the morphism (in the category of group
  schemes) $G \to G$ given by the $t$th iterate of the Frobenius.

  Now by applying the same argument of \cite[II, 3.19]{Jantzen} one
  obtains the identity
$$
\nabla((p^t-1) + (a-1)p^t) \cong \nabla(p^t-1)\otimes \nabla(a-1)^{F^t}.
$$
Noting that $\nabla(a-1)=L(a-1)$ we obtain from Steinberg's tensor
product theorem that $\nabla(ap^t-1) = L(ap^t-1).$ The equality
$\Delta(ap^t-1) = L(ap^t-1)$ follows by contravariant duality.
\end{pf}

As corollary to the above we have:
\begin{lem}\label{lem:4}
For $r$ as above, $L(r)\otimes E$ has a good
filtration, the sections of which are isomorphic with
$\nabla(r+1)$ and $\nabla(r-1)$.
\end{lem}

\begin{pf}
This follows from the Wang-Donkin-Mathieu theorem and Young's rule
since $L(r)\otimes E = \nabla(r) \otimes \nabla(1)$. One can also give
a direct proof using properties of the induction functor.
\end{pf}

\begin{lem} \label{lem:5}
For $r$ as above the modules $\nabla(r+1)$ and
$\nabla(r-1)$ are uniserial modules whose unique
composition series have the (respective) form
$$
[L(r_0), L(r_{t-1}), \dots, L(r_1)] 
$$
and 
\begin{align*}
[L(r_1), \dots, L(r_{t-1})] &\qquad \mbox{ for } a=1,\\ 
[L(r_1), \dots, L(r_{t-1}), L(r_t)] &\qquad \mbox{ for } a\ne 1
\end{align*}
where $r_0=r+1$ and $r_i=r+1-2p^{i-1}$ for $i \in \Z_{>0}$.
\end{lem}

\begin{pf}
  Uniserial (dual) Weyl modules for two-part partitions were
  classified in \cite[Proposition 2.2]{EH:uniserial}.  The result also
  follows for instance from the results in \cite{Doty:thesis}.
\end{pf}

We can now prove Theorems \ref{thm:LtensorL1} and \ref{thm:Structure}.
The first statement of Theorem \ref{thm:Structure} is contained already
in Lemma \ref{lem:2}.  Suppose that $r$ has the form
$(ap^t-1)+p^{t+1}k$ for integers $t>0$, $k\ge0$ and $a \in
\{1,\dots,p-1\}$. By Theorem \ref{thm:Tensor}, $M=L(r)\otimes E$ is
indecomposable with socle isomorphic to $L(r-1)$. Since simple modules
are contravariantly self-dual, so is $M$.  The description of the
module structure given in the statement of Theorem \ref{thm:Structure}
now follows from Lemmas \ref{lem:4} and \ref{lem:5}. This completes
the proof of the Theorems \ref{thm:LtensorL1} and \ref{thm:Structure}.
\noindent
%
%We consider the case $k=0$. 
%In this case by \cite[Theorem
%B]{BK:translation} we see that $L(r)\otimes E$ is indecomposable.
%Since simple modules are contravariantly self-dual, so is $L(r)\otimes
%E$. Moreover, it has (again by \cite[Theorem B]{BK:translation})
%simple socle (and head) isomorphic to $L(r+1)$.  The description of
%the module structure given in the statement of Theorem
%\ref{thm:Structure} now follows from Lemmas
%\ref{lem:4} and \ref{lem:5}.\\
%\newline
%\noindent
%
%
\begin{cor}
  For any $r \in \Z_{\ge0}$ of the form $r=ap^t-1$ where $a\in
  \{1,\dots,p-1\}$ and $t>0$ we have an isomorphism $T(r+1) \cong
  L(r)\otimes E$.
\end{cor}

\begin{pf}
  As we have seen above, $L(r)\otimes E$ has in this case a good
  filtration.  But simple modules are contravariantly self-dual, so
  $L(r)\otimes E$ is self-dual as well. Thus it has a
  $\Delta$-filtration.
\end{pf}

The preceding corollary considers the case $k=0$. Now suppose that
$k>0$. Set $\tilde{r} = r_1-kp^{t+1}$. Then $\tilde{r} \in \Z_{\ge0}$
and by the corollary $L(\tilde{r})\otimes E$ is a tilting module with
module structure as described. In particular, all of its composition
factors have highest weights that belong to the set $X_t$.  Thus by
Lemma \ref{lem:1} and Steinberg's tensor product theorem tensoring by
$L(k)^{F^{t+1}}$ induces an isomorphism of submodule lattices between
$L(\tilde{r})\otimes E$ and $L(r)\otimes E$.
\begin{cor}
  If $r \equiv -1 \pmod{p}$,
  the tensor product $L(r)\otimes E$ is either a tilting module or a
  `shift' of a tilting module.
\end{cor}
%
%
%This completes the proof
%of Theorem \ref{thm:Structure} in the general case.\\
%\newline
%\noindent
%
%
%Now we prove Theorem \ref{thm:LtensorL1}.  Part (b) follows from
%Theorem \ref{thm:Structure}. Part (a) is immediate from Steinberg's
%tensor product theorem.  To get part (c), consider first the case
%where $r$ belongs to the set $\{1,\dots,p-2\}$. Then (for instance
%from the strong linkage principle) we know that $L(r)=\nabla(r)$. Thus
%$L(r)\otimes E = \nabla(r)\otimes \nabla(1)$ has a good filtration. By
%Young's rule the sections in that filtration are $\nabla(r+1)$ and
%$\nabla(r-1)$. But both of these modules are simple Moreover, they
%belong to different blocks. Thus there can be no non-split extension
%between them, and part (c) is proved under our assumption on $r$.\\
%\newline
%\noindent
%
%
%Now suppose that $r$ has the form $a+pk$ for $a\in
%\{1,\dots,p-2\}$. Then we can tensor the module $L(r-kp)\otimes E$ by
%$L(k)^{F}$ and apply Lemma \ref{lem:1} to conclude that the claims
%of part (c) are valid in this case as well.  Theorem
%\ref{thm:LtensorL1} is proved.

\newcommand{\ot}{\otimes}
\parskip=6pt
\section{Examples: $p=2$}

In this section we give some examples of tensor product decompositions
in the case $p=2$. Unlike the situation for odd primes, at $p=2$ a
tensor product of simple modules is always indecomposable. 

In this and the following section, when we write $M=[r_1, \dots, r_k]$
we mean that $M$ is a uniserial module whose unique composition series
has the indicated form. We shall always identify simple factors with
their highest weight, and write = to indicate isomorphism.

For modules that are not uniserial, we will give whenever possible the
module structure diagram (in the sense of Alperin
\cite{Alperin:diagrams}), where the structure is depicted by a graph
in which the vertices correspond with the simple factors, and an edge
connects two vertices if and only if a non-split extension between those
compostion factors occurs as a subquotient of the module.

Module structure can be computed using known facts about good
filtrations, along with the known structure of dual Weyl modules (see
\cite{Doty:thesis}, \cite{Krop}), along with other basic facts which
can be found in \cite{Jantzen}.  We also make use of the fact that
$\Ext^1$ between two simple modules is at most $1$-dimensional.

$[1]\ot[1] = [0,2,0] = T(2)$. 

\bigskip%%degree 3

$[2]\ot[1] = [3] = T(3)$.

\bigskip
%%%degree 4

$[3]\ot[1] = ([1]\ot[2])\ot[1] = ([1]\ot[1])\ot[2]=[0,2,0]\ot[2] =
[2,0,4,0,2] = T(4)$.

$[2]\ot[2] = ([1]\ot[1])^F = [0,4,0]$.

\bigskip%%%degree 5

$[4]\ot[1] = [5]$.

$[3]\ot[2] = [1]\ot[2]\ot[2] = [1]\ot[0,4,0] = [1,5,1] = T(5)$.

\bigskip%%%degree 6

$[5]\ot[1] = [1]\ot[1]\ot[4] = [0,2,0]\ot[4] = [4,6,4]$.

$[4]\ot[2] = [6]$.

$[3]\ot[3] = [1]\ot[2]\ot[1]\ot[2] = [1]\ot[1]\ot[2]\ot[2] =
[0,2,0]\ot[0,4,0] = T(6)$, where $T(6)$ has structure: 
$$T(6) = 
\begin{minipage}{36mm}
  \setlength{\unitlength}{0.6mm}
\begin{picture}(60,40)(-30,0)
\thicklines
\drawline[40](3,3)(7,7)
\drawline[40](13,13)(17,17)

\drawline[40](-7,13)(-3,17)
\drawline[40](3,23)(7,27)

\drawline[40](-17,23)(-13,27)
\drawline[40](-7,33)(-3,37)

\drawline[40](17,23)(13,27)
\drawline[40](7,33)(3,37)

\drawline[40](-3,3)(-7,7)
\drawline[40](-13,13)(-17,17)

\put(0,0){\makebox(0,0){0}}
\put(-10,10){\makebox(0,0){4}}
\put(10,10){\makebox(0,0){2}}
\put(-20,20){\makebox(0,0){0}}
\put(0,20){\makebox(0,0){6}}
\put(20,20){\makebox(0,0){0}}
\put(-10,30){\makebox(0,0){2}}
\put(10,30){\makebox(0,0){4}}
\put(0,40){\makebox(0,0){0}}

\end{picture}
\end{minipage}
$$

\bigskip%%%degree 7

$[6]\ot[1] = [1] \ot [2] \ot [4] = [7] = T(7)$.

$[5]\ot[2]  = [1] \ot [2] \ot [4]= [7] = T(7)$.

$[4]\ot[3]  = [1] \ot [2] \ot [4]= [7] = T(7)$.

\bigskip%%%degree 8

$[7]\ot[1] = [1]\ot[1]\ot[2]\ot[4] = [0,2,0]\ot[2]\ot[4] =
[2,0,4,0,2]\ot[4] = [6,4,0,8,0,4,6] = T(8)$.

$[6]\ot[2] = [2]\ot[2]\ot[4] = [0,4,0]\ot[4] = [4,0,8,0,4]$.

$[5]\ot[3] = [1]\ot[1]\ot[2]\ot[4] = T(8)$.

$[4]\ot[4] = ([1]\ot[1])^{F^2} = [0,8,0]$.

\bigskip%%%degree 9

$[8]\ot[1] = [9]$.

$[7]\ot[2] = [1]\ot[2]\ot[2]\ot[4] = [1]\ot[0,4,0]\ot[4] =
[1]\ot[4,0,8,0,4] = [5,1,9,1,5] = T(9)$.

$[6]\ot[3] = [1]\ot[2]\ot[2]\ot[4] = T(9)$.

$[5]\ot[4] = [1]\ot[4]\ot[4] = [1]\ot[0,8,0] = [1,9,1]$.

\bigskip%%degree 10

$[9]\ot[1] = [1]\ot[1]\ot[8] = [0,2,0]\ot[8] = [8,10,8]$.

$[8]\ot[2] = [10]$.

$[7]\ot[3] = [1]\ot[1]\ot[2]\ot[2]\ot[4] = [4,6,4]\ot[0,4,0] = T(10)$,
where $T(10)$ has structure:
\vspace{4mm}
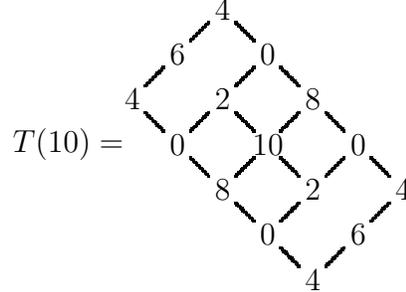
$$T(10) = 
\begin{minipage}{42mm}
  \setlength{\unitlength}{0.6mm}
\begin{picture}(60,40)(-30,0)
\thicklines
\drawline[40](3,3)(7,7)
\drawline[40](13,13)(17,17)

\drawline[40](-7,13)(-3,17)
\drawline[40](3,23)(7,27)

\drawline[40](-17,23)(-13,27)
\drawline[40](-7,33)(-3,37)

\drawline[40](17,23)(13,27)
\drawline[40](7,33)(3,37)

\drawline[40](-3,3)(-7,7)
\drawline[40](-13,13)(-17,17)

\drawline[40](-3,23)(-7,27)
\drawline[40](7,13)(3,17)

\drawline[40](13,-7)(17,-3)
\drawline[40](23,3)(27,7)
\drawline[40](23,17)(27,13)
\drawline[40](3,-3)(7,-7)

\drawline[40](-27,33)(-23,37)
\drawline[40](-17,43)(-13,47)
\drawline[40](-23,23)(-27,27)
\drawline[40](-7,47)(-3,43)

\put(0,0){\makebox(0,0){0}}
\put(-10,10){\makebox(0,0){8}}
\put(10,10){\makebox(0,0){2}}
\put(-20,20){\makebox(0,0){0}}
\put(0,20){\makebox(0,0){10}}
\put(20,20){\makebox(0,0){0}}
\put(-10,30){\makebox(0,0){2}}
\put(10,30){\makebox(0,0){8}}
\put(0,40){\makebox(0,0){0}}

\put(10,-10){\makebox(0,0){4}}
\put(20,0){\makebox(0,0){6}}
\put(30,10){\makebox(0,0){4}}

\put(-30,30){\makebox(0,0){4}}
\put(-20,40){\makebox(0,0){6}}
\put(-10,50){\makebox(0,0){4}}

\end{picture}
\end{minipage}
$$
\vspace{3mm}

$[6]\ot[4] = [2]\ot[4]\ot[4] = [2]\ot[0,8,0] = [2,10,2]$.

$[5]\ot[5] = [1]\ot[1]\ot[4]\ot[4] = [0,2,0]\ot[0,8,0] = 
\begin{minipage}{36mm}
  \setlength{\unitlength}{0.6mm}
\begin{picture}(60,40)(-30,0)
\thicklines
\drawline[40](3,3)(7,7)
\drawline[40](13,13)(17,17)

\drawline[40](-7,13)(-3,17)
\drawline[40](3,23)(7,27)

\drawline[40](-17,23)(-13,27)
\drawline[40](-7,33)(-3,37)

\drawline[40](17,23)(13,27)
\drawline[40](7,33)(3,37)

\drawline[40](-3,3)(-7,7)
\drawline[40](-13,13)(-17,17)

\drawline[40](-3,23)(-7,27)
\drawline[40](7,13)(3,17)

\put(0,0){\makebox(0,0){0}}
\put(-10,10){\makebox(0,0){8}}
\put(10,10){\makebox(0,0){2}}
\put(-20,20){\makebox(0,0){0}}
\put(0,20){\makebox(0,0){10}}
\put(20,20){\makebox(0,0){0}}
\put(-10,30){\makebox(0,0){2}}
\put(10,30){\makebox(0,0){8}}
\put(0,40){\makebox(0,0){0}}

\end{picture}
\end{minipage}
$

\bigskip%%degree 11

$[10]\ot[1] = [1]\ot[2]\ot[8] = [11]$.

$[9]\ot[2] = [1]\ot[2]\ot[8] = [11]$.

$[8]\ot[3] = [1]\ot[2]\ot[8] = [11]$.

$[7]\ot[4] = [1]\ot[2]\ot[4]\ot[4] = [3,11,3] = T(11)$.

$[6]\ot[5] = [1]\ot[2]\ot[4]\ot[4] = [3,11,3] = T(11)$.

\bigskip%%degree 12

$[11]\ot[1] = [1]\ot[1]\ot[2]\ot[8] = [2,0,4,0,2]\ot[8] = [10,8,12,8,10]$.

$[10]\ot[2] = [2]\ot[2]\ot[8] = [0,4,0]\ot[8] = [8,12,8]$.

$[9]\ot[3] = [1]\ot[1]\ot[2]\ot[8] = [10,8,12,8,10]$.

$[8]\ot[4] = [12]$.

$[7]\ot[5] = [1]\ot[1]\ot[2]\ot[4]\ot[4] = [2,0,4,0,2]\ot[0,8,0] =
T(12)$, where $T(12)$ has structure diagram:
$$T(12) = 
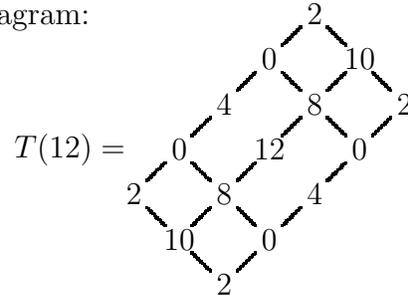
\begin{minipage}{36mm}
  \setlength{\unitlength}{0.6mm}
\begin{picture}(60,40)(-30,0)
\thicklines
\drawline[40](3,3)(7,7)
\drawline[40](13,13)(17,17)

\drawline[40](-7,13)(-3,17)
\drawline[40](3,23)(7,27)

\drawline[40](-17,23)(-13,27)
\drawline[40](-7,33)(-3,37)

\drawline[40](17,23)(13,27)
\drawline[40](7,33)(3,37)

\drawline[40](-3,3)(-7,7)
\drawline[40](-13,13)(-17,17)

\drawline[40](-3,-3)(-7,-7)
\drawline[40](-13,7)(-17,3)
\drawline[40](-23,17)(-27,13)

\drawline[40](-13,-7)(-17,-3)
\drawline[40](-23,3)(-27,7)

\drawline[40](27,27)(23,23)
\drawline[40](17,37)(13,33)
\drawline[40](7,47)(3,43)

\drawline[40](27,33)(23,37)
\drawline[40](17,43)(13,47)

\put(0,0){\makebox(0,0){0}}
\put(-10,10){\makebox(0,0){8}}
\put(10,10){\makebox(0,0){4}}
\put(-20,20){\makebox(0,0){0}}
\put(0,20){\makebox(0,0){12}}
\put(20,20){\makebox(0,0){0}}
\put(-10,30){\makebox(0,0){4}}
\put(10,30){\makebox(0,0){8}}
\put(0,40){\makebox(0,0){0}}

\put(30,30){\makebox(0,0){2}}
\put(20,40){\makebox(0,0){10}}
\put(10,50){\makebox(0,0){2}}

\put(-10,-10){\makebox(0,0){2}}
\put(-20,0){\makebox(0,0){10}}
\put(-30,10){\makebox(0,0){2}}

\end{picture}
\end{minipage}
$$
\vspace{3mm}

$[6]\ot[6] = [2]\ot[2]\ot[4]\ot[4] = [0,4,0]\ot[0,8,0] =
\begin{minipage}{36mm}
  \setlength{\unitlength}{0.6mm}
\begin{picture}(60,40)(-30,0)
\thicklines
\drawline[40](3,3)(7,7)
\drawline[40](13,13)(17,17)

\drawline[40](-7,13)(-3,17)
\drawline[40](3,23)(7,27)

\drawline[40](-17,23)(-13,27)
\drawline[40](-7,33)(-3,37)

\drawline[40](17,23)(13,27)
\drawline[40](7,33)(3,37)

\drawline[40](-3,3)(-7,7)
\drawline[40](-13,13)(-17,17)

\put(0,0){\makebox(0,0){0}}
\put(-10,10){\makebox(0,0){8}}
\put(10,10){\makebox(0,0){4}}
\put(-20,20){\makebox(0,0){0}}
\put(0,20){\makebox(0,0){12}}
\put(20,20){\makebox(0,0){0}}
\put(-10,30){\makebox(0,0){4}}
\put(10,30){\makebox(0,0){8}}
\put(0,40){\makebox(0,0){0}}

\end{picture}
\end{minipage}
$

\bigskip%%%degree13

$[12]\ot[1] = [1]\ot[4]\ot[8] = [13]$.

$[11]\ot[2] = [1]\ot[2]\ot[2]\ot[8] = [1,5,1]\ot[8] = [9,13,9]$.

$[10]\ot[3] = [1]\ot[2]\ot[2]\ot[8] = [1,5,1]\ot[8] = [9,13,9]$.

$[9]\ot[4] = [1]\ot[4]\ot[8] = [13]$.

$[8]\ot[5] = [1]\ot[4]\ot[8] = [13]$.

$[7]\ot[6] = [1]\ot[2]\ot[2]\ot[4]\ot[4] = [1,5,1]\ot[0,8,0] = T(13)$,
where $T(13)$ has structure:
$$T(13) = 
\begin{minipage}{36mm}
  \setlength{\unitlength}{0.6mm}
\begin{picture}(60,40)(-30,0)
\thicklines
\drawline[40](3,3)(7,7)
\drawline[40](13,13)(17,17)

\drawline[40](-7,13)(-3,17)
\drawline[40](3,23)(7,27)

\drawline[40](-17,23)(-13,27)
\drawline[40](-7,33)(-3,37)

\drawline[40](17,23)(13,27)
\drawline[40](7,33)(3,37)

\drawline[40](-3,3)(-7,7)
\drawline[40](-13,13)(-17,17)

\put(0,0){\makebox(0,0){1}}
\put(-10,10){\makebox(0,0){9}}
\put(10,10){\makebox(0,0){5}}
\put(-20,20){\makebox(0,0){1}}
\put(0,20){\makebox(0,0){13}}
\put(20,20){\makebox(0,0){1}}
\put(-10,30){\makebox(0,0){5}}
\put(10,30){\makebox(0,0){9}}
\put(0,40){\makebox(0,0){1}}

\end{picture}
\end{minipage}
$$

\bigskip%%%degree14

$[13]\ot[1] = [1]\ot[1]\ot[4]\ot[8]= [0,2,0]\ot[12] = [12,14,12]$

$[12]\ot[2] = [2]\ot[4]\ot[8] = [14]$.

$[11]\ot[3] = [1]\ot[1]\ot[2]\ot[2]\ot[8] = T(6)\ot[8] =
\begin{minipage}{36mm}
  \setlength{\unitlength}{0.6mm}
\begin{picture}(60,40)(-30,0)
\thicklines
\drawline[40](3,3)(7,7)
\drawline[40](13,13)(17,17)

\drawline[40](-7,13)(-3,17)
\drawline[40](3,23)(7,27)

\drawline[40](-17,23)(-13,27)
\drawline[40](-7,33)(-3,37)

\drawline[40](17,23)(13,27)
\drawline[40](7,33)(3,37)

\drawline[40](-3,3)(-7,7)
\drawline[40](-13,13)(-17,17)

\put(0,0){\makebox(0,0){8}}
\put(-10,10){\makebox(0,0){12}}
\put(10,10){\makebox(0,0){10}}
\put(-20,20){\makebox(0,0){8}}
\put(0,20){\makebox(0,0){14}}
\put(20,20){\makebox(0,0){8}}
\put(-10,30){\makebox(0,0){10}}
\put(10,30){\makebox(0,0){12}}
\put(0,40){\makebox(0,0){8}}

\end{picture}
\end{minipage}  
$

$[10]\ot[4] = [2]\ot[4]\ot[8] = [14]$.

$[9]\ot[5] = [1]\ot[1]\ot[4]\ot[8] = [12,14,12]$.

$[8]\ot[6] = [2]\ot[4]\ot[8] = [14]$.

$[7]\ot[7] = [1]\ot[1]\ot[2]\ot[2]\ot[4]\ot[4] =
[0,2,0]\ot[0,4,0]\ot[0,8,0] =T(14)$.  
%We are not able to give the complete structure diagram of $T(14)$. 
$T(14)$ has $27$ composition factors (counting multiplicities) and its
diagram is a cube consisting of $27$ vertices (standing on a vertex
labelled $0$) with the following three `layers':
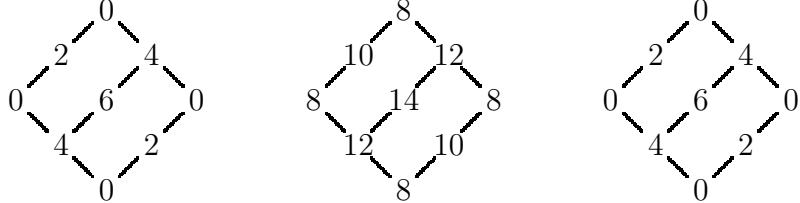
$$
\begin{array}{lll}
\begin{minipage}{36mm}
  \setlength{\unitlength}{0.6mm}
\begin{picture}(60,40)(-30,0)
\thicklines
\drawline[40](3,3)(7,7)
\drawline[40](13,13)(17,17)

\drawline[40](-7,13)(-3,17)
\drawline[40](3,23)(7,27)

\drawline[40](-17,23)(-13,27)
\drawline[40](-7,33)(-3,37)

\drawline[40](17,23)(13,27)
\drawline[40](7,33)(3,37)

\drawline[40](-3,3)(-7,7)
\drawline[40](-13,13)(-17,17)

\put(0,0){\makebox(0,0){0}}
\put(-10,10){\makebox(0,0){4}}
\put(10,10){\makebox(0,0){2}}
\put(-20,20){\makebox(0,0){0}}
\put(0,20){\makebox(0,0){6}}
\put(20,20){\makebox(0,0){0}}
\put(-10,30){\makebox(0,0){2}}
\put(10,30){\makebox(0,0){4}}
\put(0,40){\makebox(0,0){0}}

\end{picture}
\end{minipage}
&
\begin{minipage}{36mm}
  \setlength{\unitlength}{0.6mm}
\begin{picture}(60,40)(-30,0)
\thicklines
\drawline[40](3,3)(7,7)
\drawline[40](13,13)(17,17)

\drawline[40](-7,13)(-3,17)
\drawline[40](3,23)(7,27)

\drawline[40](-17,23)(-13,27)
\drawline[40](-7,33)(-3,37)

\drawline[40](17,23)(13,27)
\drawline[40](7,33)(3,37)

\drawline[40](-3,3)(-7,7)
\drawline[40](-13,13)(-17,17)

\put(0,0){\makebox(0,0){8}}
\put(-10,10){\makebox(0,0){12}}
\put(10,10){\makebox(0,0){10}}
\put(-20,20){\makebox(0,0){8}}
\put(0,20){\makebox(0,0){14}}
\put(20,20){\makebox(0,0){8}}
\put(-10,30){\makebox(0,0){10}}
\put(10,30){\makebox(0,0){12}}
\put(0,40){\makebox(0,0){8}}

\end{picture} 
\end{minipage} 
&

\begin{minipage}{36mm}
  \setlength{\unitlength}{0.6mm}
\begin{picture}(60,40)(-30,0)
\thicklines
\drawline[40](3,3)(7,7)
\drawline[40](13,13)(17,17)

\drawline[40](-7,13)(-3,17)
\drawline[40](3,23)(7,27)

\drawline[40](-17,23)(-13,27)
\drawline[40](-7,33)(-3,37)

\drawline[40](17,23)(13,27)
\drawline[40](7,33)(3,37)

\drawline[40](-3,3)(-7,7)
\drawline[40](-13,13)(-17,17)

\put(0,0){\makebox(0,0){0}}
\put(-10,10){\makebox(0,0){4}}
\put(10,10){\makebox(0,0){2}}
\put(-20,20){\makebox(0,0){0}}
\put(0,20){\makebox(0,0){6}}
\put(20,20){\makebox(0,0){0}}
\put(-10,30){\makebox(0,0){2}}
\put(10,30){\makebox(0,0){4}}
\put(0,40){\makebox(0,0){0}}

\end{picture}
\end{minipage}

\end{array}
$$
When putting the layers together to a cube there are non-split
extensions between $0,8$ and between $2, 10$ (none between $4,12$ and
$6,14$). The cube is rotated in space so as to have the two opposite
vertices labelled $0$ at the top and the bottom, and so that the socle
layers are obtained by intersecting horizontal planes through the
cube.  We therefore obtain for $T(14)$ the diagram (only the three
visible sides of the cube are shown)

$$
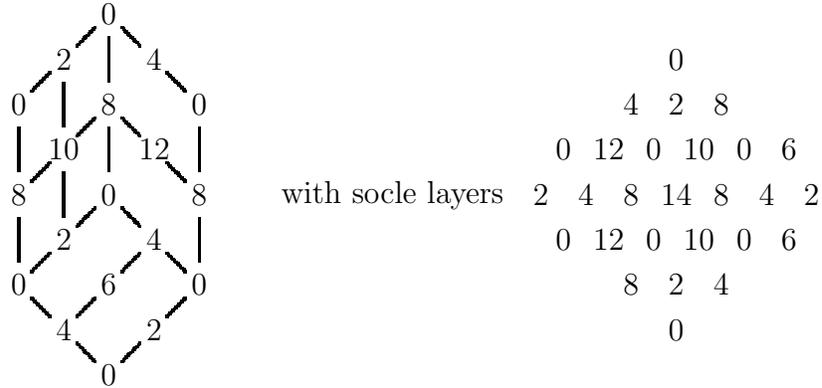
\begin{array}{lll}
\begin{minipage}{36mm}
  \setlength{\unitlength}{0.6mm}
\begin{picture}(60,80)(-30,0)
\thicklines
\drawline[40](3,3)(7,7)
\drawline[40](13,13)(17,17)

\drawline[40](-7,13)(-3,17)
\drawline[40](3,23)(7,27)

\drawline[40](-17,23)(-13,27)
\drawline[40](-7,33)(-3,37)

\drawline[40](17,23)(13,27)
\drawline[40](7,33)(3,37)

\drawline[40](-3,3)(-7,7)
\drawline[40](-13,13)(-17,17)

\drawline[40](-17,63)(-13,67)
\drawline[40](-7,73)(-3,77)

\drawline[40](-17,43)(-13,47)
\drawline[40](-7,53)(-3,57)

\drawline[40](3,77)(7,73)
\drawline[40](13,67)(17,63)
\drawline[40](3,57)(7,53)
\drawline[40](13,47)(17,43)
\drawline[40](-20,25)(-20,35)
\drawline[40](-20,45)(-20,55)
\drawline[40](-10,35)(-10,45)
\drawline[40](-10,55)(-10,65)
\drawline[40](0,45)(0,55)
\drawline[40](0,65)(0,75)
\drawline[40](20,25)(20,35)
\drawline[40](20,45)(20,55)

\put(0,0){\makebox(0,0){0}}
\put(-10,10){\makebox(0,0){4}}
\put(10,10){\makebox(0,0){2}}
\put(-20,20){\makebox(0,0){0}}
\put(0,20){\makebox(0,0){6}}
\put(20,20){\makebox(0,0){0}}
\put(-10,30){\makebox(0,0){2}}
\put(10,30){\makebox(0,0){4}}
\put(0,40){\makebox(0,0){0}}

\put(-20,40){\makebox(0,0){8}}
\put(20,40){\makebox(0,0){8}}

\put(-10,50){\makebox(0,0){10}}
\put(10,50){\makebox(0,0){12}}

\put(-20,60){\makebox(0,0){0}}
\put(0,60){\makebox(0,0){8}}
\put(20,60){\makebox(0,0){0}}
\put(-10,70){\makebox(0,0){2}}
\put(10,70){\makebox(0,0){4}}
\put(0,80){\makebox(0,0){0}}

\end{picture}
\end{minipage}
& \mbox{ with socle layers } &

\begin{minipage}{36mm}
  \setlength{\unitlength}{0.6mm}
\begin{picture}(60,60)(-30,0)
\thicklines
\put(0,0){\makebox(0,0){0}}

\put(-10,10){\makebox(0,0){8}}
\put(0,10){\makebox(0,0){2}}
\put(10,10){\makebox(0,0){4}}

\put(-25,20){\makebox(0,0){0}}
\put(-15,20){\makebox(0,0){12}}
\put(-5,20){\makebox(0,0){0}}
\put(5,20){\makebox(0,0){10}}
\put(15,20){\makebox(0,0){0}}
\put(25,20){\makebox(0,0){6}}

\put(-30,30){\makebox(0,0){2}}
\put(-20,30){\makebox(0,0){4}}
\put(-10,30){\makebox(0,0){8}}
\put(0,30){\makebox(0,0){14}}
\put(10,30){\makebox(0,0){8}}
\put(20,30){\makebox(0,0){4}}
\put(30,30){\makebox(0,0){2}}

\put(-25,40){\makebox(0,0){0}}
\put(-15,40){\makebox(0,0){12}}
\put(-5,40){\makebox(0,0){0}}
\put(5,40){\makebox(0,0){10}}
\put(15,40){\makebox(0,0){0}}
\put(25,40){\makebox(0,0){6}}

\put(-10,50){\makebox(0,0){4}}
\put(0,50){\makebox(0,0){2}}
\put(10,50){\makebox(0,0){8}}

\put(0,60){\makebox(0,0){0}}
\end{picture}
\end{minipage}
\end{array}
$$

It should be noted that $T(2^k-2)=St_{k-1} \ot St_{k-1}$ where
$St_{k-1}$ stands for the $(k-1)$th Steinberg module. These tilting
modules are projective indecomposable modules. For details
see~\cite{EH:RingelDuality}, Section 3, where it was classified when a
projective indecomposable module is a tilting module. Moreover, let
$U_i$ be defined to be the uniserial module $U_i:=L(2^{i}) \ot
L(2^{i}) = [0,2^{i+1},0]$. Then $T:=T(2^k-2)=T(2^{k-1}-2) \ot U_{k-1}$
is a $(k-1)$-dimensional hypercube with $3^{k-1}$ vertices and with
$(k-1)$ filtrations
\[
T \supset T_1 \supset T_2 \supset 0
\]
such that $T/T_1$, $T_1/T_2$ and $T_2/0$ are $(k-2)$-dimensional
hypercubes, according to tensoring $T=V_i \ot U_i$ where $V_i$ is 
given by 
\[
V_i= [1] \ot [1] \ot \ldots \ot [2^{i-1}] \ot [2^{i-1}]\ot
%  \widehat{[2^i]} \ot \widehat{[2^i]} \ot 
[2^{i+1}] \ot [2^{i+1}] \ldots \ot
  [2^{k-2}] \ot [2^{k-2}].
\]
%Here $\widehat{\mbox{\hspace{.3cm}}}$ stands for omitting the factor $[2^i]$.
%
%
In particular, in one of these filtrations, two copies of
$T(2^{k-1}-2)$ occur, one as a submodule and one as quotient module of
$T(2^k-2)$.

\section{Examples: $p=3$}

In this section we give some examples of tensor product decompositions
in the case $p=3$. Unlike the situation for $p=2$, in this case a
tensor product of simple modules is not always indecomposable.  We
follow the notational conventions introduced at the beginning of the
preceding section.  For comparing these tilting modules to projective
indecomposable modules we refer to~~\cite{EH:RingelDuality}, section
4.

\bigskip%degree2
$[1]\ot[1] = [2]\oplus[0] = T(2)\oplus T(0)$.

\bigskip%degree3
$[2]\ot[1] = [1,3,1] = T(3)$.

\bigskip%degree4
$[3]\ot[1] = [4]$.

$[2]\ot[2] = [0,4,0] \oplus [2] = T(4) \oplus T(2)$.

\bigskip%degree5
$[4]\ot[1] = [5]\oplus[3]=T(5)\oplus [3]$.

$[3]\ot[2] = [5] = T(5)$. 

\bigskip%degree6
$[5]\ot[1] = [3]\ot[2]\ot[1] = [3]\ot[1,3,1] = T(6)$, where $T(6)$
has structure:
$$T(6) = 
\begin{minipage}{12mm}
  \setlength{\unitlength}{0.6mm}
\begin{picture}(20,20)(-20,0)
\thicklines
\drawline[40](3,3)(7,7)
\drawline[40](-3,3)(-7,7)

\drawline[40](-7,13)(-3,17)
\drawline[40](7,13)(3,17)

\put(0,0){\makebox(0,0){4}}
\put(-10,10){\makebox(0,0){6}}
\put(10,10){\makebox(0,0){0}}
\put(0,20){\makebox(0,0){4}}
\end{picture}
\end{minipage}
$$

$[4]\ot[2] = [3]\ot[1]\ot[2] = T(6)$ (above). 

$[3]\ot[3] = ([1]\ot[1])^F = ([2]\oplus[0])^F = [6]\oplus [0]$.

\bigskip %degree7
$[6]\ot[1] = [7]$.

$[5]\ot[2] = [3]\ot[2]\ot[2] = [3]\ot([0,4,0]\oplus[2]) = T(7)\oplus
[5]$, where $T(7)$ has structure:
$$T(7) = 
\begin{minipage}{12mm}
  \setlength{\unitlength}{0.6mm}
\begin{picture}(20,20)(-20,0)
\thicklines
\drawline[40](3,3)(7,7)
\drawline[40](-3,3)(-7,7)

\drawline[40](-7,13)(-3,17)
\drawline[40](7,13)(3,17)

\put(0,0){\makebox(0,0){3}}
\put(-10,10){\makebox(0,0){7}}
\put(10,10){\makebox(0,0){1}}
\put(0,20){\makebox(0,0){3}}
\end{picture}
\end{minipage}
$$

$[4]\ot[3] = [3]\ot[3]\ot[1] = ([6]\oplus[0])\ot[1] = [7]\oplus[1]$.

\bigskip %degree 8

$[7]\ot[1] = [6]\ot[1]\ot[1] = [6]\ot([2]\oplus[0]) =
[8]\oplus[6]=T(8) \oplus [6]$.

$[6]\ot[2] = [8]=T(8)$.

$[5]\ot[3] = [3]\ot[3]\ot[2] = ([6]\oplus[0])\ot[2] =
[8]\oplus[2]=T(8) \oplus T(2)$.

$[4]\ot[4] = [3]\ot[3]\ot[1]\ot[1] = ([6]\oplus[0])\ot([2]\oplus[0]) 
= [8]\oplus[6]\oplus[2]\oplus[0]$.

\bigskip %degree 9

$[8]\ot[1] = [6]\ot[2]\ot[1] = [6]\ot[1,3,1] = [7,3,9,3,7] = T(9)$. 

$[7]\ot[2] = [6]\ot[1]\ot[2] = T(9)$ (above).

$[6]\ot[3] = ([2]\ot[1])^F = [1,3,1]^F = [3,9,3]$. 

$[5]\ot[4] = [3]\ot[3]\ot[2]\ot[1] = ([6]\oplus[0])\ot [1,3,1] 
= T(9) \oplus T(3)$.

\bigskip %degree 10

$[9]\ot[1] = [10]$.

$[8]\ot[2] = [6]\ot[2]\ot[2] = [6]\ot([0,4,0]\oplus[2]) 
= [6,4,10,4,6]\oplus[8] = T(10)\oplus T(8)$. 

$[7]\ot[3] = [6]\ot[3]\ot[1] = [3,9,3]\ot[1] = [4,10,4]$. 

$[6]\ot[4] = [6]\ot[3]\ot[1] = [4,10,4]$. 

$[5]\ot[5] = [3]\ot[3]\ot[2]\ot[2] = ([6]\oplus[0])\ot([0,4,0]\oplus[2])
= [6,4,10,4,6]\oplus[8] \oplus [0,4,0] \oplus [2] = T(10)\oplus T(8) 
\oplus T(4) \oplus T(2)$.

\bigskip %degree 11

$[10]\ot[1] = [9]\ot[1]\ot[1] = [9]\ot([2]\oplus[0]) = [11]\oplus[9]$.

$[9]\ot[2] = [11]$.

$[8]\ot[3] = [6]\ot[3]\ot[2] = [3,9,3]\ot[2] = [5,11,5] = T(11)$.

$[7]\ot[4] = [6]\ot[3]\ot[1]\ot[1] = [3,9,3]\ot([2]\oplus[0]) = 
T(11) \oplus   [3,9,3]$.

$[6]\ot[5] = [6]\ot[3]\ot[2] = [5,11,5] = T(11)$.

\bigskip %degree 12

$[11]\ot[1] = [9]\ot[2]\ot[1] = [9]\ot[1,3,1] = [10,12,10]$. 

$[10]\ot[2] = [9]\ot[2]\ot[1] = [10,12,10]$. 

$[9]\ot[3] = [12]$.

$[8]\ot[4] = [6]\ot[3]\ot[2]\ot[1] = [3,9,3]\ot[1,3,1] = T(12)$,
where $T(12)$ has structure: 
$$T(12) = 
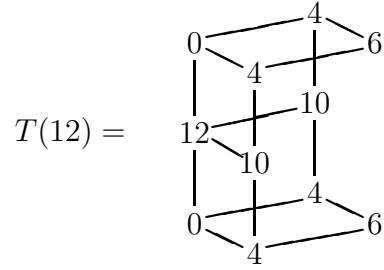
\begin{minipage}{20mm}
  \setlength{\unitlength}{0.8mm}
\begin{picture}(20,40)(-20,0)
\thicklines

\drawline[40](2,1)(18,4)
\drawline[40](-2,1)(-8,4)
\drawline[40](-8,6)(8,9)
\drawline[40](12,9)(18,6)

\drawline[40](-2,16)(-7,19)
\drawline[40](-7,21)(7,24)

\drawline[40](2,31)(18,34)
\drawline[40](-2,31)(-8,34)
\drawline[40](-8,36)(8,39)
\drawline[40](12,39)(18,36)

\drawline[40](0,3)(0,12)
\drawline[40](0,18)(0,27)
\drawline[40](-10,8)(-10,17)
\drawline[40](-10,23)(-10,32)
\drawline[40](10,13)(10,22)
\drawline[40](10,28)(10,37)

\put(0,0){\makebox(0,0){4}}
\put(0,15){\makebox(0,0){10}}
\put(-10,5){\makebox(0,0){0}}
\put(20,5){\makebox(0,0){6}}
\put(10,10){\makebox(0,0){4}}
\put(-10,20){\makebox(0,0){12}}
\put(0,30){\makebox(0,0){4}}
\put(20,35){\makebox(0,0){6}}
\put(-10,35){\makebox(0,0){0}}
\put(10,25){\makebox(0,0){10}}
\put(10,40){\makebox(0,0){4}}

\end{picture}
\end{minipage}
$$

%$$T(12) = 
%\begin{minipage}{20mm}
%  \setlength{\unitlength}{0.8mm}
%\begin{picture}(20,40)(-20,0)
%\thicklines
%\drawline[40](3,3)(7,7)
%\drawline[40](0,3)(0,7)
%\drawline[40](-3,3)(-7,7)
%\drawline[40](-10,13)(-10,17)
%\drawline[40](-7,13)(-3,17)
%\drawline[40](0,13)(0,17)
%\drawline[40](2,12)(7,17)
%\drawline[40](10,13)(10,17)
%\drawline[40](-10,23)(-10,27)
%\drawline[40](0,23)(0,27)
%\drawline[40](3,23)(7,27) 
%\drawline[40](-7,23)(-2,28)
%\drawline[40](10,23)(10,27)
%\drawline[40](7,33)(3,37)
%\drawline[40](0,33)(0,37)
%\drawline[40](-7,33)(-3,37)
%\put(0,0){\makebox(0,0){4}}
%\put(-10,10){\makebox(0,0){10}}
%\put(0,10){\makebox(0,0){0}}
%\put(10,10){\makebox(0,0){6}}
%\put(-10,20){\makebox(0,0){4}}
%\put(0,20){\makebox(0,0){12}}
%\put(10,20){\makebox(0,0){4}}
%\put(-10,30){\makebox(0,0){6}}
%\put(0,30){\makebox(0,0){0}}
%\put(10,30){\makebox(0,0){10}}
%\put(0,40){\makebox(0,0){4}}
%\end{picture}
%\end{minipage}
%$$

$[7]\ot[5] = [6]\ot[3]\ot[2]\ot[1] = T(12)$ (see above). 

$[6]\ot[6] = ([2]\ot[2])^F = [0,4,0]^F \oplus [2]^F = [0,12,0] \oplus [6]$.

\bigskip %degree 13

$[12]\ot[1] = [9]\ot[3]\ot[1] = [13]$.

$[11]\ot[2] = [9]\ot[2]\ot[2] = [9]\ot([0,4,0]\oplus [2]) = [9,13,9]
\oplus [11]$.

$[10]\ot[3] = [9]\ot[3]\ot[1] = [13]$.

$[9]\ot[4] = [9]\ot[3]\ot[1] = [13]$.

$[8]\ot[5] = [6]\ot[3]\ot[2]\ot[2] = [3,9,3]\ot([0,4,0]\oplus [2])
= [3,9,3]\ot[0,4,0] \oplus [3,9,3]\ot[2] = T(13) \oplus T(11)$, where
$T(13)$ has structure:

$$T(13) = 
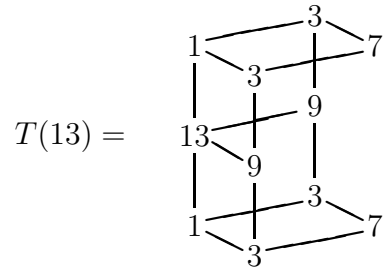
\begin{minipage}{20mm}
  \setlength{\unitlength}{0.8mm}
\begin{picture}(20,40)(-20,0)
\thicklines

\drawline[40](2,1)(18,4)
\drawline[40](-2,1)(-8,4)
\drawline[40](-8,6)(8,9)
\drawline[40](12,9)(18,6)

\drawline[40](-2,16)(-7,19)
\drawline[40](-7,21)(7,24)

\drawline[40](2,31)(18,34)
\drawline[40](-2,31)(-8,34)
\drawline[40](-8,36)(8,39)
\drawline[40](12,39)(18,36)

\drawline[40](0,3)(0,12)
\drawline[40](0,18)(0,27)
\drawline[40](-10,8)(-10,17)
\drawline[40](-10,23)(-10,32)
\drawline[40](10,13)(10,22)
\drawline[40](10,28)(10,37)

\put(0,0){\makebox(0,0){3}}
\put(0,15){\makebox(0,0){9}}
\put(-10,5){\makebox(0,0){1}}
\put(20,5){\makebox(0,0){7}}
\put(10,10){\makebox(0,0){3}}
\put(-10,20){\makebox(0,0){13}}
\put(0,30){\makebox(0,0){3}}
\put(20,35){\makebox(0,0){7}}
\put(-10,35){\makebox(0,0){1}}
\put(10,25){\makebox(0,0){9}}
\put(10,40){\makebox(0,0){3}}

\end{picture}
\end{minipage}
$$

%$$T(13) = 
%\begin{minipage}{20mm}
%  \setlength{\unitlength}{0.8mm}
%\begin{picture}(20,40)(-20,0)
%\thicklines
%\drawline[40](3,3)(7,7)
%\drawline[40](0,3)(0,7)
%\drawline[40](-3,3)(-7,7)
%\drawline[40](-10,13)(-10,17)
%\drawline[40](-7,13)(-3,17)
%\drawline[40](0,13)(0,17)
%\drawline[40](2,12)(7,17)
%\drawline[40](10,13)(10,17)
%\drawline[40](-10,23)(-10,27)
%\drawline[40](0,23)(0,27)
%\drawline[40](3,23)(7,27)
%\drawline[40](-7,23)(-2,28)
%\drawline[40](10,23)(10,27)
%\drawline[40](7,33)(3,37)
%\drawline[40](0,33)(0,37)
%\drawline[40](-7,33)(-3,37)

%\put(0,0){\makebox(0,0){3}}
%\put(-10,10){\makebox(0,0){9}}
%\put(0,10){\makebox(0,0){1}}
%\put(10,10){\makebox(0,0){7}}
%\put(-10,20){\makebox(0,0){3}}
%\put(0,20){\makebox(0,0){13}}
%\put(10,20){\makebox(0,0){3}}
%\put(-10,30){\makebox(0,0){7}}
%\put(0,30){\makebox(0,0){1}}
%\put(10,30){\makebox(0,0){9}}
%\put(0,40){\makebox(0,0){3}}
%\end{picture}
%\end{minipage}
%$$

$[7]\ot[6] = [6]\ot[6]\ot[1] = ([0,12,0] \oplus [6])\ot[1] 
= [1,13,1] \oplus [7]$.

\bigskip %degree 14

$[13]\ot[1] = [9]\ot[3]\ot[1]\ot[1] = [9]\ot[3]\ot([2]\oplus[0]) = 
[14] \oplus [12]$. 

$[12]\ot[2] = [9]\ot[3]\ot[2] = [14]$.

$[11]\ot[3] = [9]\ot[3]\ot[2] = [14]$.

$[10]\ot[4] = [9]\ot[3]\ot[1]\ot[1] = [14] \oplus [12]$. 

$[9]\ot[5] = [9]\ot[3]\ot[2] = [14]$.

$[8]\ot[6] = [6]\ot[6]\ot[2] = ([0,12,0] \oplus [6])\ot[2] 
= [2,14,2] \oplus [8] = T(14) \oplus T(8)$. 

$[7]\ot[7] = [6]\ot[6]\ot[1]\ot[1] = ([0,12,0]\oplus[6]) \ot
([2]\oplus[0]) = [2,14,2] \oplus [8] \oplus [0,12,0] \oplus [6] 
= T(14) \oplus T(8) \oplus [0,12,0] \oplus [6]$.

\bigskip %degree 15

$[14]\ot[1] = [9]\ot[3]\ot[2]\ot[1] = [9]\ot T(6) = 
\begin{minipage}{12mm}
  \setlength{\unitlength}{0.6mm}
\begin{picture}(20,20)(-20,0)
\thicklines
\drawline[40](3,3)(7,7)
\drawline[40](-3,3)(-7,7)

\drawline[40](-7,13)(-3,17)
\drawline[40](7,13)(3,17)

\put(0,0){\makebox(0,0){13}}
\put(-10,10){\makebox(0,0){15}}
\put(10,10){\makebox(0,0){9}}
\put(0,20){\makebox(0,0){13}}
\end{picture}
\end{minipage}
$

$[13]\ot[2] = [9]\ot[3]\ot[2]\ot[1]$ (see above).

$[12]\ot[3] = [9]\ot[3]\ot[3] = [9] \ot ([6]\oplus[0]) = [15]\oplus[9]$.

$[11]\ot[4] = [9]\ot[3]\ot[2]\ot[1]$ (see above).

$[10]\ot[5] = [9]\ot[3]\ot[2]\ot[1]$ (see above).

$[9]\ot[6] = ([3]\ot[2])^F = [5]^F = [15]$.

$[8]\ot[7] = [6]\ot[6]\ot[2]\ot[1] = ([0,12,0]\oplus[6])\ot[1,3,1]
= [0,12,0]\ot[1,3,1] \oplus [6]\ot[1,3,1] = T(15) \oplus T(9)$,
where $T(15)$ has structure: 

$$T(15) = 
\begin{minipage}{20mm}
  \setlength{\unitlength}{0.8mm}
\begin{picture}(20,40)(-20,0)
\thicklines

\drawline[40](-2,1)(-8,4)
\drawline[40](-8,6)(8,9)

\drawline[40](3,16)(17,19)
\drawline[40](-3,17)(-7,19)
\drawline[40](-7,21)(7,24)
\drawline[40](13,24)(17,22)

\drawline[40](-2,31)(-8,34)
\drawline[40](-8,36)(8,39)

\drawline[40](0,3)(0,12)
\drawline[40](0,18)(0,27)
\drawline[40](-10,8)(-10,17)
\drawline[40](-10,23)(-10,32)
\drawline[40](10,13)(10,22)
\drawline[40](10,28)(10,37)

\put(0,0){\makebox(0,0){1}}
\put(0,15){\makebox(0,0){13}}
\put(-10,5){\makebox(0,0){3}}
\put(10,10){\makebox(0,0){1}}
\put(-10,20){\makebox(0,0){9}}
\put(0,30){\makebox(0,0){1}}
\put(-10,35){\makebox(0,0){3}}
\put(10,25){\makebox(0,0){13}}
\put(10,40){\makebox(0,0){1}}
\put(20,20){\makebox(0,0){15}}

\end{picture}
\end{minipage}
$$

%$$T(15) = 
%\begin{minipage}{40mm}
%  \setlength{\unitlength}{0.8mm}
%\begin{picture}(40,40)(-30,0)
%\thicklines
%\drawline[40](3,3)(7,7)
%\drawline[40](-3,3)(-7,7)
%\drawline[40](-13,13)(-17,17)
%\drawline[40](13,13)(17,17)
%\drawline[40](-10,13)(-10,17)
%\drawline[40](10,13)(10,17)
%\drawline[40](-8,11)(8,19)
%\drawline[40](-17,23)(-13,27)
%\drawline[40](17,23)(13,27)
%\drawline[40](10,23)(10,27)
%\drawline[40](-8,29)(8,21)
%\drawline[40](8,29)(-8,21)
%\drawline[40](3,37)(7,33)
%\drawline[40](-3,37)(-7,33)
%
%\put(0,0){\makebox(0,0){1}}
%\put(-10,10){\makebox(0,0){13}}
%\put(10,10){\makebox(0,0){3}}
%\put(-20,20){\makebox(0,0){1}}
%\put(-10,20){\makebox(0,0){15}}
%\put(10,20){\makebox(0,0){9}}
%\put(20,20){\makebox(0,0){1}}
%\put(-10,30){\makebox(0,0){3}}
%\put(10,30){\makebox(0,0){13}}
%\put(0,40){\makebox(0,0){1}}
%\end{picture}
%\end{minipage}
%$$

\bigskip %degree 16

$[15]\ot[1] = [9]\ot[6]\ot[1] = [16]$.

$[14]\ot[2] = [9]\ot[3]\ot[2]\ot[2] = [9]\ot(T(7)\oplus[5]) = 
\begin{minipage}{12mm}
  \setlength{\unitlength}{0.6mm}
\begin{picture}(20,20)(-20,0)
\thicklines
\drawline[40](3,3)(7,7)
\drawline[40](-3,3)(-7,7)

\drawline[40](-7,13)(-3,17)
\drawline[40](7,13)(3,17)

\put(0,0){\makebox(0,0){12}}
\put(-10,10){\makebox(0,0){16}}
\put(10,10){\makebox(0,0){10}}
\put(0,20){\makebox(0,0){12}}
\end{picture}
\end{minipage}
\hspace{9mm}\oplus [14]$.

$[13]\ot[3] = [9]\ot[3]\ot[3]\ot[1] = [9]\ot([6]\oplus[0])\ot[1] 
= [16]\oplus[10]$. 

$[12]\ot[4] = [16]\oplus[10]$. 

$[11]\ot[5] = [9]\ot[3]\ot[2]\ot[2] = [9]\ot(T(7)\oplus [5]) = 
[9]\ot T(7) \oplus [14]$ (see $[14]\ot[2]$ above). 

$[10]\ot[6] = [9]\ot[6]\ot[1] = [16]$. 

$[9]\ot[7] = [9]\ot[6]\ot[1] = [16]$. 

$[8]\ot[8] = [6]\ot[6]\ot[2]\ot[2] = ([0,12,0]\oplus[6])\ot([0,4,0]\oplus[2])
= [0,12,0]\ot[0,4,0] \oplus [0,12,0]\ot[2] \oplus [6]\ot[0,4,0] \oplus [6]\ot[2]
= T(16) \oplus T(14) \oplus T(10) \oplus T(8)$, where $T(16)$ has structure:

$$T(16) = 
\begin{minipage}{20mm}
  \setlength{\unitlength}{0.8mm}
\begin{picture}(20,40)(-20,0)
\thicklines

\drawline[40](-2,1)(-8,4)
\drawline[40](-8,6)(8,9)

\drawline[40](3,16)(17,19)
\drawline[40](-3,17)(-7,19)
\drawline[40](-7,21)(7,24)
\drawline[40](13,24)(17,22)

\drawline[40](-2,31)(-8,34)
\drawline[40](-8,36)(8,39)

\drawline[40](0,3)(0,12)
\drawline[40](0,18)(0,27)
\drawline[40](-10,8)(-10,17)
\drawline[40](-10,23)(-10,32)
\drawline[40](10,13)(10,22)
\drawline[40](10,28)(10,37)

\put(0,0){\makebox(0,0){0}}
\put(0,15){\makebox(0,0){12}}
\put(-10,5){\makebox(0,0){4}}
\put(10,10){\makebox(0,0){0}}
\put(-10,20){\makebox(0,0){10}}
\put(0,30){\makebox(0,0){0}}
\put(-10,35){\makebox(0,0){4}}
\put(10,25){\makebox(0,0){12}}
\put(10,40){\makebox(0,0){0}}
\put(20,20){\makebox(0,0){16}}

\end{picture}
\end{minipage}
$$

%$$T(16) = 
%\begin{minipage}{40mm}
%  \setlength{\unitlength}{0.8mm}
%\begin{picture}(40,40)(-30,0)
%\thicklines
%\drawline[40](3,3)(7,7)
%\drawline[40](-3,3)(-7,7)
%\drawline[40](-13,13)(-17,17)
%\drawline[40](13,13)(17,17)
%\drawline[40](-10,13)(-10,17)
%\drawline[40](10,13)(10,17)
%\drawline[40](-8,11)(8,19)
%\drawline[40](-17,23)(-13,27)
%\drawline[40](17,23)(13,27)
%\drawline[40](10,23)(10,27)
%\drawline[40](-8,29)(8,21)
%\drawline[40](8,29)(-8,21)
%\drawline[40](7,29)(-17,21)
%\drawline[40](3,37)(7,33)
%\drawline[40](-3,37)(-7,33)
%\put(0,0){\makebox(0,0){0}}
%\put(-10,10){\makebox(0,0){12}}
%\put(10,10){\makebox(0,0){4}}
%\put(-20,20){\makebox(0,0){0}}
%\put(-10,20){\makebox(0,0){16}}
%\put(10,20){\makebox(0,0){10}}
%\put(20,20){\makebox(0,0){0}}
%\put(-10,30){\makebox(0,0){4}}
%\put(10,30){\makebox(0,0){12}}
%\put(0,40){\makebox(0,0){0}}
%\end{picture}
%\end{minipage}
%$$

%\newpage

\small\parskip=0pt
\par Mathematical and Computer Sciences
\par Loyola University Chicago
\par Chicago, Illinois 60626 U.S.A.
\medskip 
\par E-mail: {\tt doty@math.luc.edu} 

\bigskip 
\par Mathematics and Computer Science
\par University of Leicester
\par University Road
\par Leicester  LE1 7RH U.K.

\medskip
\par E-mail: {\tt A.Henke@mcs.le.ac.uk} 

\end{document}